\documentclass[11pt, twoside]{article}
\usepackage{latexsym}
\usepackage{amsmath}
\usepackage{amssymb}
\usepackage[all]{xy}
\usepackage{amsfonts}
\usepackage{verbatim}
\usepackage{amsthm}
\usepackage{mathrsfs}
\usepackage{epsfig}
\usepackage{xy}
\usepackage{array}
\usepackage{stmaryrd}
\usepackage{graphicx,color}
\usepackage{xcolor}
\usepackage[colorlinks=true,linkcolor=blue,citecolor=blue]{hyperref}
\usepackage{tikz}
\usetikzlibrary{arrows,calc}
\usepackage{etex}
\usepackage{mathdots}
\usepackage{float}
\usepackage{graphics}
\usepackage{pdflscape}
\usepackage{bm}
\usepackage{anysize,hyperref}
\input xypic
\xyoption{all}

\usepackage[perpage,symbol]{footmisc}
\usepackage{setspace}
\def\dr{\ar@{->}[r]}

\SelectTips{cm}{10}

\begin{document}
\title{\Large{\bf Classifying Nichols algebras over classical Weyl groups$^\bigstar$\footnotetext{\hspace{-1em}$^\bigstar$Weicai Wu was supported by Hunan Provincial Natural Science Foundation of China (Grant No: 2020JJ5210),Key Laboratory of Mathematical and Statistical Model(Guangxi Normal University), Education Department of Guangxi Zhuang Autonomous Region and Panyue Zhou was supported by the National Natural Science Foundation of China (Grant No. 11901190).}}}
\medskip
\author{Weicai Wu and Panyue Zhou}

\date{}

\maketitle
\def\blue{\color{blue}}
\def\red{\color{red}}

\newtheorem{theorem}{Theorem}[section]
\newtheorem{lemma}[theorem]{Lemma}
\newtheorem{corollary}[theorem]{Corollary}
\newtheorem{proposition}[theorem]{Proposition}
\newtheorem{conjecture}{Conjecture}
\theoremstyle{definition}
\newtheorem{definition}[theorem]{Definition}
\newtheorem{question}[theorem]{Question}
\newtheorem{remark}[theorem]{Remark}
\newtheorem{remark*}[]{Remark}
\newtheorem{example}[theorem]{Example}
\newtheorem{example*}[]{Example}
\newtheorem{condition}[theorem]{Condition}
\newtheorem{condition*}[]{Condition}
\newtheorem{construction}[theorem]{Construction}
\newtheorem{construction*}[]{Construction}

\newtheorem{assumption}[theorem]{Assumption}
\newtheorem{assumption*}[]{Assumption}

\baselineskip=17pt
\parindent=0.5cm
\vspace{-6mm}

\begin{abstract}
\baselineskip=16pt
In this article, we show that conjugacy classes of classical Weyl groups $W(B_{n})$ and $W(D_{n})$ are of \textit{type  D}. Consequently, we obtain that Nichols algebras of irreducible Yetter-Drinfeld modules over the classical Weyl groups $\mathbb W_{n}$ ($n\geq5$) are infinite dimensional.\\[0.2cm]
\end{abstract}

\pagestyle{myheadings}
\markboth{\rightline {\scriptsize   W. Wu and P. Zhou}}
         {\leftline{\scriptsize Classifying Nichols algebras over classical Weyl groups}}

\section{Introduction}
The theory of Nichols algebras is motivated by the Hopf algebra theory.
It has interesting applications to other research fields such as Kac-Moody Lie superalgebras, and conformal field theory and quantum groups. In particular, the field of quantum groups and pointed Hopf
algebras over finite abelian groups has seen great progress due to the
work of Andruskiewitsch, Schneider, Heckenberger, see \cite {4},\cite {7},\cite {8}, etc. In any area of mathematics the classification of all objects is very important. In Hopf algebra theory, the classification of all finite dimensional pointed Hopf algebras is a tough question \cite{0}.
In \cite{5'}, they classify finite-dimensional pointed Hopf algebras with abelian coradical, and show that they are cocycle deformations of the associated graded Hopf algebra, moreover, the complete result on the non-abelian groups has not yet been obtained. Given a group $ G$, an important step to classify all finite-dimensional pointed Hopf algebras $H$ with group-like $G(H) = G$ is to determine all the pairs $(\mathcal O, \rho )$ such that the associated Nichols algebra $  \mathfrak B(\mathcal O, \rho )$ is finite-dimensional. If the Nichols algebra is infinite dimensional, then the Hopf algebra is so, it is often useful to discard those pairs such that $\dim  \mathfrak B (\mathcal O, \rho ) = \infty$. There are properties of the conjugacy class $\mathcal O$ which imply $\dim  \mathfrak B (\mathcal O, \rho )  = \infty$ for any $\rho $, one of which is the property of being of \textit{type  D}.

Nichols algebras of braided vector spaces  $(\mathbb CX,  c_{q})$,  where $X$ is a
rack and $q$ is a 2-cocycle in $X$, were studied in \cite {3}, among others. It was shown \cite {1},\cite {2},\cite {5} that Nichols algebras  $\mathfrak B({\mathcal O}_
\sigma, \rho)$ over symmetry groups $\mathbb S_{n}$ (isomorphism with $W(A_{n-1})$) are infinite dimensional,  except in three cases; some open cases were subsequently treated by Fantino in \cite{5''}.
In \cite{10},\cite{11},\cite{12} they showed that except in several cases conjugacy classes of classical Weyl groups $W(B_{n})$
(isomorphism with $W(C_{n})$) and $W(D_{n})$ are of \textit{type  D},
as a consequence, they show
that Nichols algebras of a large list of irreducible Yetter-Drinfeld ({\rm YD} in short) modules over
those Weyl groups are infinite dimensional.
However, the classification has not been completed for Nichols algebras over
classic Weyl groups $W(B_{n})$ and $W(D_{n})$.

We will study of classify finite-dimensional Nichols algebra over the classical Weyl groups. Note that $\mathbb Z_{2}^{n} \rtimes \mathbb S_{n}$
is isomorphic to Weyl groups $W(B_{n})$ and $W(C_{n})$ of $B_{n}$ and
$C_{n}$ for $n>2$. If $K_{n}= \{a \in
     \mathbb Z_{2}^{n} \mid \ a = ({a_{1}}, {a_{2}}, \cdots, {a_{n}}) $  $\hbox {
with   } a_{1} +a_{2} + \cdots + a_{n} \equiv 0  \ ({\rm mod } \ 2) \}$,
then $K_{n} \rtimes \mathbb S_n$ is isomorphic to Weyl group
$W(D_{n})$ of $D_{n}$ for $n>3$ and $ K_{n} \rtimes  \mathbb{S}_{n}$  is a subgroup of
$\mathbb Z_{2}^{n} \rtimes \mathbb S_{n}$. Let $\mathbb W_{n}$ denote $H_{n}\rtimes  \mathbb{S}_{n}$ where
$H_{n}:=  K_{n}$  or $\mathbb Z_{2}^{n}$.

For any $a\in H_{n}$ with $a = ({a_{1}}, {a_{2}},\cdots, {a_{n}})$ and $\sigma, \tau \in \mathbb{S}_{n}$, $\sigma (a) :=( {a_{\sigma^{-1}(1)}}, {a_{\sigma^{-1}(2)}}, \cdots,{a_{\sigma^{-1}(n)}} )$,

$(a, \sigma )^{- 1} = (- (a _{\sigma (1)}, a_{\sigma (2)}, \cdots  , a_{\sigma (n)}), \sigma^{ - 1}) = (- \sigma ^{ - 1}(a), \sigma ^{- 1})$,

$(b, \tau  )(a, \sigma )(b, \tau  )^{- 1} = (b + \tau  (a) -  \tau  \sigma \tau ^{ - 1}(b), \tau  \sigma \tau ^{ - 1})$.

Sometimes we abbreviate $a \sigma$ as $(a, \sigma)$. The other notations are the same as in \cite {11}.

This article is organized as follows. In Section 2 we provide some preliminaries.
In Section 3  we examine the rack and conjugacy classes of $\mathbb W_{n}$ and prove that Nichols algebras of irreducible {\rm YD} modules over the classical Weyl groups $\mathbb W_{n}$ ($n\geq5$) are infinite dimensional.

\section{Preliminaries}\label {s1}
Let $\mathbb{K}$ be an algebraically closed field of characteristic zero, $G$ a finite group, and $V$ a finite dimensional vector space over  $\mathbb{K} $. For $s\in G$ and  $(\rho,  V) \in  \widehat {G^s}$,  here is a
precise description of the {\rm YD} module $M({\mathcal O}_s,
\rho)$,  introduced in \cite {6}. Let $t_1 = s, t_2,   \cdots,
t_{m}$ be a numeration of ${\mathcal O}_s$,  which is a conjugacy
class containing $s$,   and let $g_i\in G$ such that $g_i \rhd s :=
g_i s g_i^{-1} = t_i$ for all $1\le i \le m$. Then $M({\mathcal
O}_s,  \rho) = \oplus_{1\le i \le m}g_i\otimes V$. Let $g_iv :=
g_i\otimes v \in M({\mathcal O}_s, \rho)$,  $1\le i \le m$,  $v\in V$.
If $v\in V$ and $1\le i \le m $,  then the action of $h\in G$ and the
coaction are given by
\begin {eqnarray} \label {e0.11}
\delta(g_iv) = t_i\otimes g_iv,  \qquad h\cdot (g_iv) =
g_j(\nu _i(h)\cdot v),  \end {eqnarray}
 where $hg_i = g_j\nu _i(h)$,  for
unique  $1\le j \le m$ and $\nu _i(h)\in G^s$. The explicit formula for
the braiding is then given by
\begin{equation} \label{yd-braiding}
C(g_iv\otimes g_jw) = t_i\cdot(g_jw)\otimes g_iv =
g_{j'}(\nu _j (t_i)\cdot w)\otimes g_iv\end{equation} for any $1\le i, j\le
m$,  $v, w\in V$,  where $t_ig_j = g_{j'}\nu _j(t_i)$ for unique $j'$,  $1\le
j' \le m$ and $\nu _j(t_i) \in G^s$. Let $\mathfrak{B} ({\mathcal O}_s,
\rho )$ denote $\mathfrak{B} (M ({\mathcal O}_s,  \rho ))$.
$M({\mathcal O}_s,  \rho )$ is a simple {\rm YD} module (see \cite
 {5}).

\textit{A rack} is a pair $(X, \rhd )$,  where $ X\neq\emptyset,\rhd :  X \times X \rightarrow  X$ is an operation such that
$x\rhd (y \rhd z) = (x\rhd y)\rhd (x\rhd z)$ for any $x, y, z
\in  X$ and the map $\phi_{x}: y\rightarrow x\rhd y$
is invertible, for any $x, y\in X.$
For example, $({\mathcal O}_s^G, \rhd )$ is a rack  with $x\rhd y:= x y x^{-1}.$

Assume $R$ and $S$ are two subracks of $X$  with $R\cup S = X$, $R\cap S =\emptyset$. For any $x\in R, y\in S$, if $x\rhd y \in S$, $y\rhd x \in R$,  then $R\cup S$ is called a \textit{decomposition of subracks} of $X$. Moreover, if there exist $\sigma\in R$, $\tau\in S$ such that ${\rm sq} (\sigma, \tau) := \sigma\rhd (\tau \rhd (\sigma \rhd \tau)) \not= \tau$, then $X$ is called to be of \textit{type  D}. Notice that if $X$ is a subrack of a rack $Y$, and $X$ is of \textit{type  D}, then $Y$ is also of  \textit{type  D}.
If $G$ is a finite group and  $\mathcal{O}_{\sigma}^G$ is of \textit{type  D}, then  ${\rm dim}~\mathfrak B(\mathcal{O}_{\sigma}^G,  \rho) = \infty $ for any $\rho \in \widehat { G^\sigma}$ (see \cite[Theorem 3.6]{1}).
See \cite{3},\cite{1} for more details.

\section{Main results}\label {s2}

In this section we  prove the conjugacy classes of classical Weyl groups $W(B_{n})$ and $W(D_{n})$ are of \textit{type  D} and show that Nichols algebras of irreducible {\rm YD} modules over the classical Weyl groups $\mathbb W_{n}$ ($n\geq5$) are infinite dimensional.

We begin by stating the following result of \cite{11} and \cite{12}, whose proof
uses Theorem 8.6 of \cite{9} and Proposition 3.4 of \cite{1}.

\begin{theorem}\label{1.002}{\rm \cite[Theorem 5.4]{11}} Assume $n\ge 5$. Let $\sigma\in\mathbb S_{n}$ be of type $(1^{m_{1}},2^{m_{2}},\dots,n^{m_{n}})$ and $a\in H_{n}$ with $a\sigma\in \mathbb W_{n}$ and $\sigma\not=1$. If {\rm dim} $\mathfrak B(\mathcal{O}_{a\sigma}^{\mathbb W_{n}}, \rho)<\infty $,  then the type of $\sigma$ belongs to one in the following list:
\renewcommand{\theenumi}{\roman{enumi}}   \renewcommand{\labelenumi}{(\theenumi)}
\begin{enumerate}
\item[\rm (i)] $(2,  3)$; $(2^3);$
\item[\rm (ii)]  $(2^4);$ $(1,  2^2), (1^2, 3), (1^2, 2^2) ;$
\item[\rm (iii)]     $(1^{n-2},  2)$ and  $(1^{n-3},  3)$ $(n >5)$ with $a_i = a_j$ when $\tau (i) = i $ and $\tau (j)=j$.
\end{enumerate}
\end{theorem}

\begin{theorem}\label {1.003}{\rm \cite[Theorem 4.1]{12}}  Assume $n\ge 5$. Let $\sigma\in\mathbb S_{n}$ be of type $(1^{m_{1}},2^{m_{2}},\dots,n^{m_{n}})$ and $a\in H_{n}$ with $a\sigma\in \mathbb W_{n}$ and $\sigma\not=1$. If {\rm dim} $\mathfrak B(\mathcal{O}_{a\sigma}^{\mathbb W_{n}}, \rho)<\infty $,  then the type of $\sigma$ belongs to one in the following list:
\renewcommand{\theenumi}{\roman{enumi}}   \renewcommand{\labelenumi}{(\theenumi)}
\begin{enumerate}
\item[\rm (i)] $(2, 3)$;
\item[\rm (ii)] $(1^{2}, 3)$;
\item[\rm (iii)] $(1^{n-2},  2)$ for $n >5$ with $a_{i} = a_{j}$ when $\sigma(i) = i $ and $\sigma(j)=j$.
\end{enumerate}
\end{theorem}

\begin{lemma}\label {1.001} {\rm \cite[Lemma 3.2]{12}} Let $a=(a_{1},\ldots,a_{n}),b=(b_{1},\ldots,b_{n})$ and
$\sum\limits_{i=1}^{n}a_{i}\equiv\sum\limits_{i=1}^{n}b_{i}(mod\ 2)$. Assume
$a\sigma=(a^{(1)}\sigma_{1})(a^{(2)}\sigma_{2})\cdots(a^{(r)}\sigma_{r})$
and
$ b\tau=(b^{(1)}\tau_{1})(b^{(2)}\tau_{2})\cdots(b^{(r)}\tau_{r})$
are independent sign cycle decomposition {\rm (see {\rm \cite[Theorem 4.1]{10}})} of $a\sigma, b\tau\in \mathbb W_n$, respectively.
For all $1\leq i\leq r$, if $\sigma_{i}$ and $\tau_{i}$ have the same length, then $a\sigma$ and $b\tau$ are conjugate in $\mathbb W_{n}$.
\end{lemma}

\begin{proposition}\label {1.01} Assume  $n >3$. If $\sigma$  is of type $(1^{n-2},  2)$,  then $\mathcal O_{a\sigma} ^{\mathbb W_{n}}$ is of \textit{type  D} for all $a \in H_{n}$.
\end{proposition}

\proof Let $A:=\{(12),(13),(23)\}\subseteq S_{3}$.

{\rm (i)} Assume $\mathbb W_{n}:=K_{n} \rtimes \mathbb S_{n}$, then $\sum\limits_{i=1}^{n}a_{i}\equiv0(mod\ 2)$ for $\forall\  a\in \mathbb Z_{2}^{n}$.

If $n$ is even. Let $R:= \{((0,0,\cdots,0,0),\mu)\mid \mu\in A\}$ and $S:= \{((1,1,\cdots,1,1),\mu)\mid \mu\in A\}$. If $n$ is odd. Let $R:= \{((0,0,\cdots,0,0),\mu)\mid \mu\in A\}$ and $S:= \{((1,1,\cdots,1,0),\mu)\mid \mu\in A\}$.

{\rm (ii)} Assume $\mathbb W_{n}:=\mathbb Z_{2}^{n} \rtimes \mathbb S_{n}$.
We only need to consider $\forall\  a\in \mathbb Z_{2}^{n}$ with $\sum\limits_{i=1}^{n}a_{i}\equiv1(mod\ 2)$ by {\rm (i)}.

If $n$ is even.
Let $R:= \{((0,0,\cdots,0,1),\mu)\mid \mu\in A\}$ and $S:= \{((1,1,\cdots,1,0),\mu)\mid \mu\in A\}$. If $n$ is odd.
Let $R:= \{((0,0,\cdots,0,1),\mu)\mid \mu\in A\}$ and $S:= \{((1,1,\cdots,1,1),\mu)\mid \mu\in A\}$.

For $\forall\ a\sigma\in R,b\tau\in S$, $a\sigma$ and $b\tau$ are conjugate by Lemma \ref {1.001}. It is clear that
$R \cup S$ is a subrack decomposition of   $\mathcal O_{a\sigma}^{\mathbb W_{n}}$,
Take $\sigma=(12),\tau=(13)$, then ${\rm sq} (a\sigma, b\tau)\neq b\tau$ since
${\rm sq} (\sigma, \tau)=\sigma\neq \tau$, which implies that $R\cup S$ is of \textit{type  D}. \qed

\begin{lemma}\label {1.21} Let $G= \mathbb W_{n}$ and $a\sigma, b\tau\in G$.

{\rm (i) } {\rm (See \cite[Lemma 3.1(i)]{11})} Then  ${\rm sq} (a\sigma, b\tau)=(a+\sigma(b)+\sigma\tau(a)+\sigma\tau\sigma(b)-
\sigma\tau\sigma\tau\sigma^{-1} (a)
-\sigma\tau\sigma\tau \sigma^{-1}\tau^{-1}(b)
-\sigma\tau\sigma\tau \sigma^{-1}\tau^{-1}\sigma^{-1}(a),\sigma\tau\sigma\tau \sigma^{-1}\tau^{-1}\sigma^{-1})$.

{\rm (ii)} If $\sigma\in\mathbb S_{n}$  is of type $(1^{n-3},3)$, then
${\rm sq }(a\sigma, b\sigma) =b\sigma$ if and only if
\begin {eqnarray}\label {e11}  \sigma(a)+\sigma^{2}(a)=\sigma(b)+\sigma^{2}(b).\end {eqnarray}

{\rm (iii)} If $\sigma=\tau\mu\in\mathbb S_{5}$  is of type $(2,3)$ and $\tau^{2}=1,\mu^{3}=1$, then
${\rm sq }(a\sigma, b\sigma) =b\sigma$ if and only if
\begin {eqnarray}\label {e22}  a+\tau(a)+\tau\mu(a)+\mu^{2}(a)=b+\tau(b)+\tau\mu(b)+\mu^{2}(b).\end {eqnarray}
\end{lemma}

\proof  It is straightforward to verify. \qed

\begin{proposition}\label {1.22} If $\sigma$  is of type $(1^{2},3)$,  then
$\mathcal O_{a\sigma}^{\mathbb W_{5}}$ is of \textit{type  D} for all $a \in H_{5}$.
\end{proposition}

\proof  Let $\sigma=(123)$ without lost generality. We have that
the $1$-st,  $2$-nd, $3$-rd components of (\ref {e11}) are
\begin {eqnarray}\label {e33}  (a_2+a_3,  a_1+a_3,  a_1+a_2) =  (b_2+b_3,  b_1+b_3,  b_1+b_2). \end {eqnarray}
If $a =(1,1,0,0,0)$ and $b=(1,0,0,0,1)$,  then (\ref {e33}) does not hold, if $a =(1,1,0,0,1)$ and $b=(1,0,0,0,0)$,  then (\ref {e33}) does not hold.

{\rm (i)} Assume $\mathbb W_{5}:=K_{5} \rtimes \mathbb S_{5}$, then $\sum\limits_{i=1}^{5}a_{i}\equiv0(mod\ 2)$ for $\forall\  a\in \mathbb Z_{2}^{5}$.

Let $R:= \{((0,0,0,0,0),\sigma),((1,1,0,0,0),\sigma),((1,0,1,0,0),\sigma),
((0,1,1,0,0),\sigma)\}$ and $S$
$:= \{((1,0,0,0,1),\sigma),((0,1,0,0,1),\sigma),((0,0,1,0,1),\sigma),
((1,1,1,0,1),\sigma)\}$.

{\rm (ii)} Assume $\mathbb W_{5}:=\mathbb Z_{2}^{5} \rtimes \mathbb S_{5}$.
We only need to consider $\forall\  a\in \mathbb Z_{2}^{5}$ with $\sum\limits_{i=1}^{5}a_{i}\equiv1(mod\ 2)$ by {\rm (i)}.

Let $R:= \{((0,0,0,0,1),\sigma),((1,1,0,0,1),\sigma),((1,0,1,0,1),\sigma),
((0,1,1,0,1),\sigma)\}$ and $S$
$:= \{((1,0,0,0,0),\sigma),((0,1,0,0,0),\sigma),((0,0,1,0,0),\sigma),
((1,1,1,0,0),\sigma)\}$.

For $\forall\ a\sigma\in R,b\sigma\in S$, $a\sigma$ and $b\sigma$ are conjugate by Lemma \ref {1.001}.  Let $(c,\sigma)$ denote $(a,\sigma)\rhd (b,\sigma)$, thus
$((c_{1},c_{2},c_{3},c_{4},c_{5}),\sigma)=
((a_{1}+a_{3}+b_{3},a_{1}+a_{2}+b_{1},a_{2}+a_{3}+b_{2},b_{4},b_{5}),\sigma)$.
It is clear $\sum\limits_{i=1}^{3}c_{i}\equiv\sum\limits_{i=1}^{3}b_{i}(mod\ 2)$, then
$R \cup S$ is a subrack decomposition of   $\mathcal O_{a\sigma}^{\mathbb W_{5}}$.
Consequently, $R\cup S$ is of \textit{type  D}. \qed

\begin{proposition}\label {1.23} If $\sigma$  is of type $(1,3)$,  then
$\mathcal O_{a\sigma}^{\mathbb W_{4}}$ is of \textit{type  D} for all $a \in H_{4}$.
\end{proposition}

\proof Similarly to Proposition \ref {1.22}. {\rm (i)} If $\sum\limits_{i=1}^{4}a_{i}\equiv0(mod\ 2)$ for $\forall\  a\in \mathbb Z_{2}^{4}$.
Let $R:= \{((0,0,0,0),\sigma),((1,1,0,0),\sigma),((1,0,1,0),\sigma),
((0,1,1,0),\sigma)\}$ and $S$
$:= \{((1,0,0,1),\sigma),$

\noindent $((0,1,0,1),\sigma),((0,0,1,1),\sigma),
((1,1,1,1),\sigma)\}$.
Take $a =(1,1,0,0)$ and $b=(1,0,0,1)$.
{\rm (ii)} If $\forall\  a\in \mathbb Z_{2}^{4}$ with $\sum\limits_{i=1}^{4}a_{i}\equiv1(mod\ 2)$.
Let $R:= \{((0,0,0,1),\sigma),((1,1,0,1),\sigma),((1,0,1,1),\sigma),$

\noindent $
((0,1,1,1),\sigma)\}$ and $S$
$:= \{((1,0,0,0),\sigma),((0,1,0,0),\sigma),((0,0,1,0),\sigma),
((1,1,1,0),\sigma)\}$.
Take $a =(1,1,0,1)$ and $b=(1,0,0,0)$. \qed

\begin{proposition}\label {1.23''} If $\sigma$  is of type $(4)$,  then
$\mathcal O_{a\sigma}^{\mathbb W_{4}}$ is of \textit{type  D} for all $a \in K_{4}$.
\end{proposition}

\proof Let $A:=\{(1234),(1243),(1324),(1342),(1423),(1432)\}\subseteq S_{4}$.
Let $R:= \{((0,0,0,0),\mu)\mid \mu\in A\}$ and $S:= \{((1,1,1,1),\mu)\mid \mu\in A\}$.
For $\forall\ a\sigma\in R,b\tau\in S$, $a\sigma$ and $b\tau$ are conjugate by Lemma \ref {1.001}. It is clear that
$R \cup S$ is a subrack decomposition of   $\mathcal O_{a\sigma}^{\mathbb W_{n}}$,
Take $\sigma=(1234),\tau=(1243)$, then ${\rm sq} (a\sigma, b\tau)\neq b\tau$ since
${\rm sq} (\sigma, \tau)=\sigma\neq \tau$, which implies that $R\cup S$ is of \textit{type  D}. \qed

\begin{proposition}\label {1.24} If $\sigma$  is of type $(2,3)$,  then
$\mathcal O_{a\sigma}^{\mathbb W_{5}}$ is of \textit{type  D} for all $a \in H_{5}$.
\end{proposition}

\proof Let $\sigma = (1\ 2\ 3)(4\ 5)$ without lost generality, We have that
the $1$-st,  $2$-nd, $3$-rd components of (\ref {e22}) are
\begin {eqnarray}\label {e44}  (a_2+a_3,  a_1+a_3,  a_1+a_2) =  (b_2+b_3,  b_1+b_3,  b_1+b_2). \end {eqnarray}
If $a =(1,1,0,0,0)$ and $b=(1,0,0,0,1)$,  then (\ref {e44}) does not hold, if $a =(1,1,0,0,1)$ and $b=(1,0,0,0,0)$,  then (\ref {e44}) does not hold.

{\rm (i)} Assume $\mathbb W_{5}:=K_{5} \rtimes \mathbb S_{5}$, then $\sum\limits_{i=1}^{5}a_{i}\equiv0(mod\ 2)$ for $\forall\  a\in \mathbb Z_{2}^{5}$.

Let $R:= \{((0,0,0,0,0),\sigma),((1,1,0,0,0),\sigma),((1,0,1,0,0),\sigma),
((0,1,1,0,0),\sigma),$

\noindent $((0,0,0,1,1),\sigma),((1,1,0,1,1),\sigma),((1,0,1,1,1),\sigma),
((0,1,1,1,1),\sigma)\}$

and $S:= \{((1,0,0,0,1),\sigma),((0,1,0,0,1),\sigma),((0,0,1,0,1),\sigma),
((1,1,1,0,1),\sigma),$

\noindent $((1,0,0,1,0),\sigma),((0,1,0,1,0),\sigma),((0,0,1,1,0),\sigma),
((1,1,1,1,0),\sigma)\}$

{\rm (ii)} Assume $\mathbb W_{5}:=\mathbb Z_{2}^{5} \rtimes \mathbb S_{5}$.
We only need to consider $\forall\  a\in \mathbb Z_{2}^{5}$ with $\sum\limits_{i=1}^{5}a_{i}\equiv1(mod\ 2)$ by {\rm (i)}.

Let $R:= \{((0,0,0,0,1),\sigma),((1,1,0,0,1),\sigma),((1,0,1,0,1),\sigma),
((0,1,1,0,1),\sigma)$

\noindent $((0,0,0,1,0),\sigma),((1,1,0,1,0),\sigma),((1,0,1,1,0),\sigma),
((0,1,1,1,0),\sigma)\}$

and $S:= \{((1,0,0,0,0),\sigma),((0,1,0,0,0),\sigma),((0,0,1,0,0),\sigma),
((1,1,1,0,0),\sigma),$

\noindent $((1,0,0,1,1),\sigma),((0,1,0,1,1),\sigma),((0,0,1,1,1),\sigma),
((1,1,1,1,1),\sigma)\}$.

For $\forall\ a\sigma\in R,b\sigma\in S$, $a\sigma$ and $b\sigma$ are conjugate by Lemma \ref {1.001}.  Let $(c,\sigma)$ denote $(a,\sigma)\rhd (b,\sigma)$, thus
$((c_{1},c_{2},c_{3},c_{4},c_{5}),\sigma)=
((a_{1}+a_{3}+b_{3},a_{1}+a_{2}+b_{1},a_{2}+a_{3}+b_{2},
a_{4}+a_{5}+b_{5},a_{4}+a_{5}+b_{4}),\sigma)$.
It is clear $\sum\limits_{i=1}^{3}c_{i}\equiv\sum\limits_{i=1}^{3}b_{i}(mod\ 2)$, $c_{4}+c_{5}\equiv b_{4}+b_{5}(mod\ 2)$, then
$R \cup S$ is a subrack decomposition of $\mathcal O_{a\sigma}^{\mathbb W_{5}}$.
Consequently, $R\cup S$ is of \textit{type  D}. \qed

\begin{theorem}\label{1.8}
Let $n\geq 5$. Let $\sigma\in\mathbb S_n$ be of type $(1^{m_{1}},2^{m_{2}},\dots,n^{m_{n}})$ and $a\in H_{n}$ with $a\sigma\in W_{n}$ and $\sigma\not=1$. Then $\mathcal{O}_{a\sigma}^{ \mathbb W_{n}}$ is of
\textit{type  D}.
\end{theorem}

\proof It follows from Proposition \ref {1.01}, Proposition \ref {1.22}, Proposition \ref {1.24} and Theorem \ref {1.003}. \qed
\medskip

If Nichols algebras of irreducible {\rm YD} modules over a finite group $G$ have infinite dimension, then the same holds for any finite-dimensional {\rm YD} module (see \cite{1}).

We have the following assertion by Theorem \ref {1.8}.

\begin{theorem}\label{2.6} Assume $n\ge 5$. Let $\sigma\in\mathbb S_{n}$ be of type $(1^{m_{1}},2^{m_{2}},\dots,n^{m_{n}})$ and $a\in H_{n}$ with $a\sigma\in \mathbb W_{n}$ and $\sigma\not=1$. Then {\rm dim} $\mathfrak B(\mathcal{O}_{a\sigma}^{\mathbb W_{n}}, \rho)=\infty $.
\end{theorem}

\begin{remark}
Theorem \ref{2.6} gives an affirmative answer of the conjecture due to Wu \cite[Conjecture 4.1]{12}.
\end{remark}

%
%

\textbf{Weicai Wu }\\
School of Mathematics and Statistics, Guangxi Normal University,
541004 Guilin,   Guangxi, People's Republic of China. \\
Email: weicaiwu@hnu.edu.cn\\[0.3cm]
\textbf{Panyue Zhou}\\
College of Mathematics, Hunan Institute of Science and Technology, 414006 Yueyang, Hunan,  People's Republic of China.\\
E-mail: panyuezhou@163.com

\end{document}